\documentclass[conference]{IEEEtran}

\IEEEoverridecommandlockouts
\usepackage{cite}
\usepackage{amsmath,amssymb,amsfonts}
\usepackage{mlspconf}
\usepackage{graphicx}
\usepackage{textcomp}
\usepackage{xcolor}
\usepackage{comment}
\usepackage{enumitem}
\usepackage[ruled,vlined]{algorithm2e}
\def\BibTeX{{\rm B\kern-.05em{\sc i\kern-.025em b}\kern-.08em
    T\kern-.1667em\lower.7ex\hbox{E}\kern-.125emX}}
\include{alphabet}  

\newtheorem{theorem}{Theorem}

\newtheorem{assumption}{Assumption}
\newtheorem{proposition}{Proposition}
\newtheorem{definition}{Definition}

\newcommand{\review}[1]{#1}

\copyrightnotice{979-8-3503-2411-2/23/\$31.00 {\copyright}2023 European Union}

\copyrightnotice{979-8-3503-2411-2/23/\$31.00 {\copyright}2023 IEEE}

\toappear{2023 IEEE International Workshop on Machine Learning for Signal Processing, Sept.\ 17--20, 2023, Rome, Italy}

\title{A new non-convex framework to improve asymptotical knowledge on generic stochastic gradient descent} 

\name{Jean-Baptiste Fest$^{(1)}$, Audrey Repetti$^{(2)}$, and \'Emilie Chouzenoux$^{(1)}$\thanks{This work is funded by the European Research Council Starting Grant MAJORIS ERC-2019-STG-850925.}}
\address{
(1) CVN, CentraleSupélec, Inria, Université Paris-Saclay, 9 rue Joliot Curie, Gif-sur-Yvette, France\\
(2) School of Engineering \& Physical Sciences, School of Mathematical \& Computer Schiences, \\Heriot-Watt University, EH14 4AS Edinburgh, UK
}

\begin{document}

\maketitle


\begin{abstract}
Stochastic gradient optimization methods are broadly used to minimize non-convex smooth objective functions, for instance when training deep neural networks. However, theoretical guarantees on the asymptotic behaviour of these methods remain scarce. Especially, ensuring almost-sure convergence of the iterates to a stationary point is quite challenging.
In this work, we introduce a new Kurdyka-\L ojasiewicz theoretical framework to analyze asymptotic behavior of stochastic gradient descent (SGD) schemes when minimizing non-convex smooth objectives. In particular, our framework provides new almost-sure convergence results, on iterates generated by any SGD method satisfying mild conditional descent conditions. 
We illustrate the proposed framework by means of several toy simulation examples. We illustrate the role of the considered theoretical assumptions, and investigate how SGD iterates are impacted whether these  assumptions are either fully or partially satisfied.
\end{abstract}
\smallbreak 


\section{Introduction}

We consider the unconstrained optimization problem
\begin{equation} \label{eq:pb}
\underset{\xb\in \Rbb^N}{\text{minimize}}\quad F(\xb),
\end{equation}
where $F\colon \Rbb^N \to \eR$ is a continuously differentiable function ($N\geq 1$), that is not necessarily assumed to be convex. We focus on the challenging situation when the evaluation of the gradient of $F$ is subject to (stochastic) errors, during the iterative resolution procedure. This typically arises in important scenarios of supervised machine learning, when $F$ is an expectation loss to be minimized through the access of data samples \cite{duchi2011adaptive}. In such context, the non-convexity of $F$ results from the use of nonlinear models, such as deep neural networks, for mapping the data features \cite{castera2022second,kingma2014adam}. 

The most popular approach to solve \eqref{eq:pb} in this context is probably the stochastic gradient descent (SGD)~\cite{robbins1951stochastic} and its numerous variants \cite{duchi2011adaptive,tieleman2012lecture,kingma2014adam}. SGD generates a stochastic sequence $(\xv_k)_{k\in \Nbb}$ defined as
\begin{equation}\label{eq:SGD_scheme}
    (\forall k \in \Nbb)\quad \xv_{k+1} = \xv_k - \alpha_k \gv_k,
\end{equation}
where $(\gv_k)_{k\in \Nbb}$ is typically a random process defined on a probability space $\left(\Omega, \Fc, \Pbb \right)$ aiming at approximating the true gradient $\left(\nabla F(\xv_k)\right)_{k\in \Nbb}$. Moreover, $(\alpha_k)_{k\in \Nbb}$ corresponds to a positive stepsize sequence. Practical applications of such SGD schemes to supervised learning can be found for instance in~ \cite{bottou2018optimization,khaled2020better}.

In general, SGD schemes satisfy conditional descent properties on the sequence $(F(\xv_k))_{k\in \Nbb}$, with respect to the natural filtration $(\Fc_k)_{k \in  \Nbb} = \left( \sigma( \xv_0, \dots, \xv_k) \right)_{k\in \Nbb}$. Such properties can be obtained under technical assumptions on the noise on the approximated gradients and on the sequence $(\alpha_k)_{k\in \Nbb}$~\cite{robbins1971convergence,gadat2017stochastic,bertsekas2000gradient}. 
In particular, for most SGD schemes, $(F(\xv_k))_{k\in \Nbb}$ satisfies an almost-supermartingale condition, and converges to a finite limit \cite{robbins1971convergence}. In addition, the gradient process $(\nabla F(\xv_k))_{k\in \Nbb}$ (or a sub-sequence of it) converges to zero \cite{bertsekas2000gradient}.

However, usually little is known on the asymptotic behaviour of the generated process $(\xv_k)_{k \in \Nbb}$ itself. The machine learning literature typically focuses instead on the search for (fast) convergence rates~\cite{roux2012stochastic, johnson2013accelerating, defazio2014saga, defazio2014finito}, and is often limited to the convex (or even strictly convex) setting. 
The lack of study of stochastic algorithms in a non-convex framework in the literature calls for developing new theory in this context. In particular, it is necessary to leverage existing deterministic approaches, typically relying on Kurdyka-\L ojasiewicz (KL) properties \cite{kurdyka1998gradients,bolte2014proximal}. KL properties have initially been introduced to solve gradient flow problems in a continuous setting \cite{bolte2007lojasiewicz,absil2005convergence}. Nevertheless, it has proven to be particularly efficient to improve convergences guarantees of discrete optimization schemes when no convexity assumptions are made \cite{attouch2009convergence, bolte2014proximal, chouzenoux2014variable, chouzenoux2016block, repetti2021}.

In this work, we introduce a new theoretical framework to prove almost sure (a.s.) convergence of the stochastic process $(\xv_k)_{k\in \Nbb}$ generated by SGD schemes of the form of \eqref{eq:SGD_scheme}, to a critical point of $F$, when $F$ is not necessarily convex. We show that this result applies to any SGD schemes holding mild conditional descent properties, when $F$ satisfies a Kurdyka-\L ojasiewicz (KL) property \cite{kurdyka1998gradients,bolte2014proximal}. 
We further empirically investigate how the considered assumptions practically control the convergence behaviour of SGD schemes, on some toy simulation examples.
The current paper relies on our recent preprint \cite{chouzenoux2023kurdyka}. The originality of the current work is to specialize the study to SGD, and to present a comprehensive numerical illustration for the results.

The remainder of the paper is organized as follow. Section \ref{sec2} introduces the theoretical tools required. We present our main theoretical contribution in Section \ref{sec3}. Numerical experiments are presented in Section \ref{sec4}, and Section \ref{sec5} concludes the work.  

\section{Theoretical Background}
\label{sec2}

\subsection{KL property}

One key mathematical tool for proving convergence of non-convex deterministic optimization schemes is the KL property. It has been initially introduced by \L ojasiewicz \cite{lojasiewicz1963propriete} and  Kurdika \cite{kurdyka1998gradients}, and has been at the core of major methodological developments in non-convex optimization analysis, in the last decades, starting with the seminal papers from~\cite{bolte2007lojasiewicz,attouch2010proximal}. 

The definition of KL property is given below.

\begin{definition}\label{def_KL}\emph{(KL property)}
A differentiable function $F : \Rbb^{N} \to \Rbb$  satisfies the \emph{Kurdyka-\L ojasiewicz} (KL) property on $E \subset\Hc$, if for every $\widetilde{\xb} \in E$, there exists a neighbourhood $V$ of $\widetilde{\xb}$, $\zeta > 0$ and $\varphi \in \Phi_{\zeta}$ such that $$\|\nabla F(\xb)\|\varphi'\left(F(\xb)-F(\widetilde{\xb})\right) \ge 1,$$ for every $\xb \in  V $ satisfying $$ 0 < F(\xb)-F(\widetilde{\xb}) < \zeta.$$ 
\end{definition}

\subsection{Convergence analysis under KL property}

Definition \ref{def_KL} was initially motivated in a continuous setting, through the gradient flow analysis \cite{absil2005convergence}. Indeed,  KL property promotes finite gradient trajectories converging to the origin. As a discrete counterpart of gradient flows, gradient descent (GD) algorithms are thus expected to follow a similar behaviour. 

Let us consider Problem \eqref{eq:pb}, where $F$ satisfies KL property. Let us build the sequence $(\xb_k)_{k\in \Nbb}$ generated by a (deterministic) GD algorithm of the form of
\begin{equation}
    (\forall k \in \Nbb)\quad
    \xb_{k+1} = \xb_k - \alpha_k \nabla F(\xb_k),
\end{equation}
where, for every $k\in \Nbb$, $\alpha_k>0$ is a stepsize.
If $\left(\nabla F(\xb_k)\right)_{k\in \Nbb}$ and $\left(\xb_{k+1} - \xb_k\right)_{k\in \Nbb}$ are proportional, and if 
$\sum_{k=0}^{+\infty} \|\nabla F(\xb_k) \|<+\infty$ (e.g., if $\alpha_k$ is small enough and $F$ is Lipschitz smooth),
then we can deduce from KL property that  
$$\sum_{k=0}^{+\infty} \|\xb_{k+1} - \xb_k\|<+\infty.$$ 
This result can then be used to deduce that $(\xb_k)_{k\in \Nbb}$ is a convergent Cauchy sequence. Detailed examples of convergence analysis of gradient-based schemes under KL property can be found in \cite{attouch2009convergence}.

In the non-convex case, KL property thus allows to show that the limit point of $(\xb_k)_{k\in \Nbb}$ exists and cancels the gradient (i.e., stationary point), under mild requirements such as descent conditions, on $\left(F(\xb_k)\right)_{k\in \Nbb}$. Note that local convergence results to global solutions can also be obtained, when a good initialization (i.e., close enough to a global minimum of $F$) is considered \cite{jin2017escape}. 

\smallbreak

\subsection{Uniformized KL property}

Definition \ref{def_KL} might be sometimes limited, as it is too `local' to be easily manipulated. Recently, \cite{bolte2014proximal} introduced an alternative version of the KL property that we introduce in the following Theorem.

\begin{definition}
\label{th_uniform_KL} \emph{(Uniformized KL property)}
Let $C$ be a compact set in $\Hc$ and $F \colon \Hc \to \Rbb$ be a differentiable function constant on $C$, satisfying the KL property on $C$. Then, there exist $(\varepsilon, \zeta) \in (0,+\infty)^2$ and $\varphi \in \Phi_{\zeta}$ such that 
\begin{equation}\label{eq:ineqKL_uniform}
    (\forall \overline{\xb} \in C)(\forall \xb\in \Hc) \quad
    \|\nabla F(\xb)\| \, \varphi'\left(F(\xb)-F(\overline{\xb})\right) \geq 1,
\end{equation}
when $d(\xb,C)<\varepsilon$ and $0<F(\xb)-F(\overline{\xb}) <\zeta$.
\end{definition}

\smallbreak

A typical usage of Definition~\ref{th_uniform_KL} is when $\left(F(\xb_k)\right)_{k\in \Nbb}$ is decreasing and $F$ coercive. Then,  $\left(F(\xb_k)\right)_{k\in \Nbb}$ actually converges to a finite limit while $(\xb_k)_{k\in \Nbb}$ is guaranteed to be bounded. As a consequence, taking $C$ in Definition~\ref{th_uniform_KL} as the set of cluster points of $(\xb_k)_{k\in \Nbb}$ allows \eqref{eq:ineqKL_uniform} to be verified for any iterate starting from a certain rank.

\section{Convergence of SGD for non-convex objectives}
\label{sec3}

\subsection{Generic SGD scheme}

Let us consider SGD schemes of the form of \eqref{eq:SGD_scheme} for solving \eqref{eq:pb}. We denote $(\xv_k)_{k\in \Nbb}$, a stochastic process generated by~\eqref{eq:SGD_scheme}. 
In our study, we assume that the process satisfies the two following conditions. 

\begin{assumption}
\label{ass:H1}
$F$ is coercive and $\beta$-Lipschitz differentiable on $\mathbb{R}^N$.
\end{assumption}

\begin{assumption}
\label{ass:H3}
    $F$ satisfies the KL property on the set of critical points of $F$. Moreover, this set can be written as the finite union of non-empty disjoint compacts subset. 
\end{assumption}

\smallbreak

Assumption~\ref{ass:H1} is a classical hypothesis usually made in the field of differentiable optimization context \cite{nocedal1999numerical}. In particular, it ensures the existence of a minimal value of $F$, denoted by $F_{\min}$.
On the contrary, Assumption~\ref{ass:H3} is specific to our non-convex context, as we do not have any convexity assumption on $F$.  
Omitting some technical details here, Assumption~\ref{ass:H3} is essential for us to obtain convergence results directly on the iterates following a similar strategy as those conducted, e.g., in \cite{attouch2009convergence, bolte2014proximal, chouzenoux2014variable}, but generalized to our stochastic framework. 
Note also that the geometry imposed on the critical set is only slightly constraining.

\subsection{Gradient approximation assumptions}

Before establishing our main convergence result, we first introduce technical conditions on the stochastic approximations $(\gv_k)_{k \in \Nbb}$ of the gradients $(\nabla F(\xv_k))_{k \in \Nbb}$, involved in the SGD updates.

\begin{assumption}\label{ass:H2}
There exists three deterministic non-negative sequences $(a_k)_{k\in \Nbb}, (b_k)_{k\in \Nbb},(c_k)_{k\in \Nbb}$ such that 
\begin{enumerate}[label=\roman*)]
\item \label{ass:H2_1} $\sum_{k=0}^{+\infty}\alpha_k \left(\sqrt{a_k}+\sqrt{c_k}\right)< +\infty$, 
\item \label{ass:H2_2} $\inf_{k\in \Nbb} \alpha_k\left(1-\frac{\alpha_kb_k\beta}{2}\right)>0$, 
\smallbreak
and, for every $k\in \Nbb$,
\medbreak
\item \label{ass:H2_3} $0<\sqrt{\frac{b_{k+1}}{b_k}} \frac{2-\alpha_k b_k\beta}{2-\alpha_{k+1} b_{k+1}\beta}\left(1+\frac{\alpha_{k+1}^2a_{k+1}\beta}{2} \right)\leq 1 \notag$ 
\item \label{ass:H2_4} $\Ebb_k[\gv_k] = \nabla F(\xv_k)$
\item \label{ass:H2_5} $\Ebb_k[\|\gv_k\|^2] \leq a_k (F(\xv_k)-F_{\min})+b_k\|\nabla F(\xv_k)\|^2 + c_k \notag$. 
\end{enumerate} 
\end{assumption}

\smallbreak

Although Assumption \ref{ass:H2} may seem quite demanding, it actually gathers several typical scenarios:
\begin{itemize}
\item Assumption \ref{ass:H2}-\ref{ass:H2_4} and \ref{ass:H2}-\ref{ass:H2_5} are relative to the two first moments of the noise on the gradient term. Assumption \ref{ass:H2}-\ref{ass:H2_4} classically requires the noise to be unbiased, which Assumption \ref{ass:H2}-\ref{ass:H2_5} is a mild condition on the noise variance, generalizing many encountered in the literature~\cite{schmidt2013fast,bottou2018optimization,vaswani2019fast, khaled2020better}. A non-zero $(a_k)_{k\in \Nbb}$ sequence typically models cases when $F$ has a gradient confusion bound \cite{sankararaman2020impact}. 
\item 
Assumptions \ref{ass:H2}-\ref{ass:H2_1} and \ref{ass:H2}-\ref{ass:H2_2} are classical summability and non-zero rules, controling the terms in Assumption \ref{ass:H2}-\ref{ass:H2_5}. Assumption \ref{ass:H2}-\ref{ass:H2_2} assumes a non-vanishing stepsize in the SGD update, and requires $(b_k)_{k\in \Nbb}$ to be bounded.  
\item Assumption \ref{ass:H2}-\ref{ass:H2_3} is a non-increasing condition which naturally appears in our convergence proof. 
\end{itemize}

\smallbreak

As shown in \cite{chouzenoux2023kurdyka}, Assumption~\ref{ass:H2} guarantees that the process satisfies some conditional descent properties. As a consequence, $(F(\xv_k)_{k\in \Nbb})_{k\in \Nbb}$ converges a.s. to an a.s finite random variable $\FD_{\infty}$, and that $\left(\nabla F(\xv_k)\right)_{k\in \Nbb}$ almost surely converges to $\zerob_N$ (i.e., the zero vector of $\mathbb{R}^N$).

\subsection{Proposed KL analysis for stochastic framework}

In the context of stochastic optimization, the use of KL property is challenging. A first idea would be to apply the uniformized KL property (Definition \ref{th_uniform_KL}) to any trajectory $(\xv_k(\omega))_{k\in \Nbb}$, for every $\omega \in \Omega$. 
However, by doing so, $(\epsilon, \zeta)$ and $ \varphi$ would be random objects whose analysis is very delicate, dependent on $\omega$. 
For instance, conditional expectation operations would become tricky, and measurability of $(\Fc_k)_{k\in \Nbb}$ would not be straightforward. 

\smallbreak

To overcome this challenge, we proposed in \cite{chouzenoux2023kurdyka} a new extension of Theorem~\ref{th_uniform_KL}, better adapted to a stochastic optimization framework.    

\begin{proposition}\label{prop_stoKL}
Under Assumptions \ref{ass:H1}-\ref{ass:H2}, there exists a bounded concave function $\varphi$ and an a.s. finite positive discrete random variable $\KD$ such that $$\|\nabla F(\xv_k)\|  \varphi'(F(\xv_k) - \FD_{\infty})\geq 1$$ a.s. for every $k > \KD$.
\end{proposition}

The advantage of Proposition \ref{prop_stoKL} (whose detailed proof is given in \cite{chouzenoux2023kurdyka}), lies in the random variable $\KD$ which concentrates all the stochastic information. As such, this new tool tends to overcome some technical obstacles raised in KL-based convergence analysis \cite{chouzenoux2023kurdyka}, and allows to build a new convergence theorem, that we present hereafter in the SGD case.

\subsection{Convergence result}

We now introduce our main contribution, which is the almost-sure convergence result for the generic SGD scheme~\eqref{eq:SGD_scheme}.
Let us denote by $\FD_{\infty}$ the almost-sure finite limit of $\left( F(\xv_k) \right)_{k\in \Nbb}$ and $\Ebb_k[.]$ the conditional expectation operator relative to $\Fc_k$, for $k\in \Nbb$ \review{(i.e., for a given integrable or positive random variable, $\Ebb_k[.]$ corresponds to its best approximation regarding all information available on the process from state $0$ to state $k$). We introduce, for all $(k,\gamma) \in \Nbb \times (0,1)$, the event:}

\vspace*{-0.5cm}

\begin{multline}\label{set_Xi}
    \Xi_{\gamma,k}:= \bigg\{\FD_{\infty}< F(\xv_k) \text{ and } \\[-0.2cm] 
    \Big|(\FD_{\infty}-F_{\min})-\frac{2}{2+\alpha_k^2 a_k\beta}\Ebb_k[\review{F(\xv_{k+1})}-F_{\min}] +\frac{\alpha_k^2c_k\beta}{2} \Big| \\[-0.2cm]
    \leq \gamma\big(1-\frac{\alpha_k b_k \beta}{2}\big) \|\nabla F(\xv_k)\|^2
    \bigg\}.
\end{multline}
\smallbreak

 \begin{theorem}\label{thKL}
Under Assumptions \ref{ass:H1}, \ref{ass:H3} and \ref{ass:H2}, if there exists $\gamma_0 \in (0,1)$ such that $\Pbb(\Xi_{\gamma_0,k})=1$ for all $k\in \Nbb$, then $(\xv_k)_{k\in \Nbb}$ almost-surely converges to a critical point $\xv_\infty$ of~$F$.
\end{theorem}

\smallbreak

\review{Equality $\Pbb(\Xi_{\gamma_0,k})=1$ (for all $k\in \Nbb$) supposes that process $(\xv_k)_{k\in \Nbb}$ is well-built enough to verify a suitable descent condition and to approach its limit $F_{\infty}$ from above. Moreover, it also reflects a predictability condition; the conditional decreasing as well as the difference $\FD_{\infty}-F_{\min}$ shall be controlled with respect to the evolution of the gradient norm. }

The complete proof for Theorem~\ref{thKL} can be found in \cite{chouzenoux2023kurdyka}. It relies on the use of Proposition~\ref{prop_stoKL} as a cornerstone to establish the finite length of $(\xv_k)_{k\in \Nbb}$ almost surely. 
Up to our knowledge, Theorem \ref{thKL} is one of the first results ensuring the almost convergence of a stochastic gradient type iterates in a non-convex setting.

\smallbreak

Table~\ref{tab_schemes} gives a few examples of state-of-the-art schemes directly verifying our specific Assumption~\ref{ass:H2}. This table is not exhaustive, and Assumption~\ref{ass:H2} could be verified by other algorithms, e.g., \cite{castera2022second, chouzenoux2022sabrina} (see also \cite{chouzenoux2023kurdyka} for proximal algorithms).   

\begin{table}[!t]
    \centering\footnotesize
    \begin{tabular}
{|c|c|c|c|c|c|}
    \hline
    Scheme  & $\alpha_k$ & $a_k$ & $b_k$ & $c_k$   &  Ass. 3-\ref{ass:H2_1} and  \ref{ass:H2_2} ? \\  
    \hline
     SGD \cite{schmidt2013fast} & $\alpha$ & $0$ & $B^2$ & $0$ & Yes  \\ 
     \hline
     SGD \cite{combettes2016stochastic} & $\lambda_k\gamma_k$ &  $0$ & $1+\tau_k$ & $\zeta_k$ & $\inf\limits_{k\in \Nbb}\gamma_k \lambda_k >0$  \\ 
    \hline
     BFGS \cite{wang2017stochastic}  &$\alpha_k$& $0$ & $1$ & $\sigma^2 m_k^{-1}$ & $\sum\limits_{k=0}^{+\infty} m_k^{-1} < +\infty$ \\ \hline
    BFGS \cite{meng2020stochastic} 
    & $\eta$ & $0$ & $\rho$ & $0$ & Yes \\
    \hline 
    \end{tabular}
    \caption{Examples of state-of-the-art algorithms satisfying conditions of Theorem~\ref{thKL}, hence ensuring their convergence in a non-convex setting. For the sake of readability, we perused the same notations as the authors in their articles.}
   \label{tab_schemes}
\end{table}

\section{Numerical illustrations}
\label{sec4}

In this section we conduct some experiments on a non-convex scalar problem, so as to illustrate the behaviour of SGD algorithm when Assumptions~\ref{ass:H1}, \ref{ass:H3} and \ref{ass:H2} are satisfied by function $F$, and its moments approximation gradient sequence. To do so, we proceed by gradually increasing the complexity of noise structure. 

\smallbreak

\begin{figure}[!t]
    \centering
\includegraphics[height = 4.5cm, width=9cm]{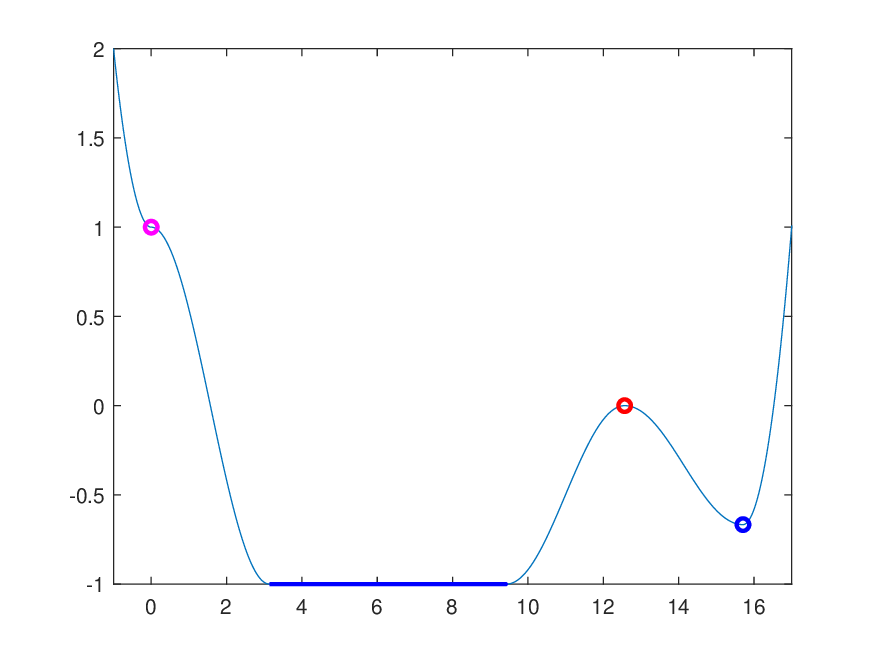}

    \vspace*{-0.5cm}

    \caption{Graph of function $F$ in \eqref{eq:fct_F_1D}. 
    Global minimizers correspond to $[\pi, 3\pi]$. 
    Point $x=5\pi$ in blue (resp.  $x=4\pi$ in in red) corresponds to a local minimizer (resp. maximizer), and $x=0$ is a saddle point.}
    \label{fig:graphF}
\end{figure}

All along this study, we work with the following $F$ function. 
\begin{equation}\label{eq:fct_F_1D}
    F : x \in \Rbb \mapsto 
    \begin{cases}
    x^2+1 \text{ if } x<0,\\
    \cos(x) \text{ if } 0\leq x<\pi, \\
    -1  \text{ if } \pi\leq x<3\pi, \\
    1/2 \cos(x) - 1/2  \text{ if } 3\pi\leq x<4\pi, \\
    1/3 \cos(x) - 1/3  \text{ if } 4\pi\leq x<5\pi, \\
    x^2 - 10\pi x +25\pi^2 -2/3 \text{ if } x\geq 5\pi. 
    \end{cases}
\end{equation}
The graph of function $F$ in \eqref{eq:fct_F_1D} is illustrated in Figure~\ref{fig:graphF}. This function is non-convex, but Lipschitz-differentiable with Lipschitz constant equals to $\beta=2$. Hence Assumption~\ref{ass:H1} holds. Moreover, since the graph of $F$ is semi-algebraic, Assumption \ref{ass:H3} is also verified \cite{bolte2007lojasiewicz}. 

In our simulations, we first numerically verify that $(\xD_k)_{k\in \Nbb}$ converges almost surely to a stationary point of $F$, denoted by $\xD_\infty$.
Then we analyse more specifically how $\xD_\infty$ is approximated by $(\xD_k)_{k\in \Nbb}$, ideally in such a way that there exists $\gamma_0 \in (0,1)$ such that $\Pbb(\Xi_{\gamma_0,k})=1$ for all $k\in \Nbb$.     
In particular, we run experiments considering different settings for the sequences $(a_k)_{k\in \Nbb}$, $(b_k)_{k\in \Nbb}$ and $(c_k)_{k\in \Nbb}$ appearing in Assumption~\ref{ass:H2}.

\subsection{Experiments under $a_k:=0,c_k:=0$}
\label{sec:3:A}

In this section we investigate the case where the only non-zero sequence (except for the stepsize) interfering in Assumption~\ref{ass:H2} is $(b_k)_{k\in \Nbb}$. 
Such kind of noise assumption is generally considered as a baseline in the literature  of stochastic optimization, as it is typically verified for the usual constant stepsize SGD scheme when $F$ satisfies the Strong Growth condition \cite{schmidt2013fast}. 
Assumptions~\ref{ass:H2} \ref{ass:H2_1}-\ref{ass:H2_3} are then verified as soon as $(b_k)_{k\in \Nbb}$ is a non-increasing sequence. 
As the latter have to be bounded to fulfill Assumption~\ref{ass:H2}-\ref{ass:H2_2}, it thus becomes equivalent to take $b_k \equiv b$, for every $k \in \Nbb$, to also verify Assumption~\ref{ass:H2}-\ref{ass:H2_5}. 
\smallbreak

The approximation sequence is generated empirically so as to satisfy both Assumption~\ref{ass:H2}-\ref{ass:H2_4} and \ref{ass:H2_5}. Specifically, we set, for every $k\in \Nbb$, $\gD_k = \eD_k^1 \nabla F(\xD_k)$, where $\eD_k^1$ is sampled uniformly in $\left[1-\sqrt{3b-1}, 1+\sqrt{3b-1}\right]$. 
For the stepsize, we set, for every $k\in \Nbb$, $\alpha_k = \alpha$ such that $\alpha b = \beta^{-1}$. 
Finally, we introduce two different perturbation levels to test the SGD algorithm in particularly extreme cases: a moderate perturbation $b=10$, and an excessively high perturbation $b=10^3$.   

\begin{table}[!t]
    \centering
    \begin{tabular}{||c||c||c|c|}
    \hline 
         $x_0$ & $b$ & Nature of $\xD_{\infty}$ & $F(\xD_k)>\FD_{\infty}~\forall k\in \Nbb$ ? 
        \\ \hline 
        $-1/2$ & $10$  & Saddle point & Yes if $F(\xD_k)\neq F(\xD_{\infty})$\\
         & $10^3$  & Saddle point  & Yes if $F(\xD_k)\neq F(\xD_{\infty})$ \\
         \hline 
          $1$ & $10$  & Global Min & Yes if $F(\xD_k)\neq F(\xD_{\infty})$ \\
          & $10^3$  & Global Min & Yes  \\
         \hline
          $4\pi + \epsilon $ & $10$  & Global or Local Min & Yes if $F(\xD_k)\neq F(\xD_{\infty})$\\
          & $10^3$  & Local Max & No  \\
         \hline 
    \end{tabular}
    \caption{Asymptotical behaviours of $(\xD_k)_{k\in \Nbb}$ under $\Ebb_k\left[\gD_k^2\right] \leq b  F'(\xD_k)^2~(k\in \Nbb)$ }
    \label{tab1}
\end{table}
Table~\ref{tab1} shows the behaviour of process $(\xD_k)_{k\in \Nbb}$ considering three different initializations $x_0 \in \{-1/2, 1, 4\pi+\epsilon\}$. The first initialization $x_0=-1/2$ is located on the left of saddle point $x=0$ (see graph of $F$ in Figure~\ref{fig:graphF}). The second initialization $x_0=1$ is between the saddle point and the interval $[\pi, 3\pi]$ of global minimizers. And the last one $x_0=4\pi+\epsilon$ is in a small neighborhood of the local maximizer $x=4\pi$, taking $\epsilon = 10^{-5}$. 
In most scenarios, $(\xD_k)_{k\in \Nbb}$ converges to a stationary point so that $F(\xD_k)_{k\in \Nbb}$ remains above its limit. The only tricky case arises for high perturbation level when $x_0=4\pi + \epsilon$, which seems to be too close to local maximizer $x=4\pi$. 

\subsection{Experiments under  $a_k:=0$}
\label{sec:3:B}

In this section, we no longer impose $(c_k)_{k\in \Nbb}$ to be equal to zero. Such a situation is regularly encountered as a first relaxed version of the noise resulting from the Strong Growth condition \cite{schmidt2013fast}. One typical situation is when, for $k\in \Nbb$, the difference between $\gD_k$ and $F'(\xD_k)$ follows a Gaussian distribution that remains independent from $\Fc_k$. 
We adopt such model to conduct our investigation. More specifically, we choose, for all $k\in \Nbb$, $\gD_k = \eD_k^1  F'(\xD_k) + \eD_k^2$. Here, $\eD_k^2$ is normally distributed, with zero-mean and standard deviation $\sigma_k >0$, and does not depend on process $(\xD_k)_{k\in \Nbb}$ so as to have $c_k = 2\sigma_k^2$. 
Moreover, $\eD_k^1$ keeps the same properties as in Section~\ref{sec:3:A}. 
In order to easily verify Assumption \ref{ass:H2} \ref{ass:H2_1}-\ref{ass:H2_3} , we set constant sequences $(\alpha_k)_{k\in \Nbb}$ and $ (b_k)_{k\in \Nbb}$, equal to $\alpha$ and $b$, respectively. In practice, we choose $b=10$ with $\alpha b = \beta^{-1}$, and $\sigma_k := \sigma / (k+1)^{(1+\varepsilon)}$ for $\sigma\in \{10 ,10^2\}$.    

\begin{figure}[!t]
    \centering
    \includegraphics[height = 4.5cm, width=9cm]{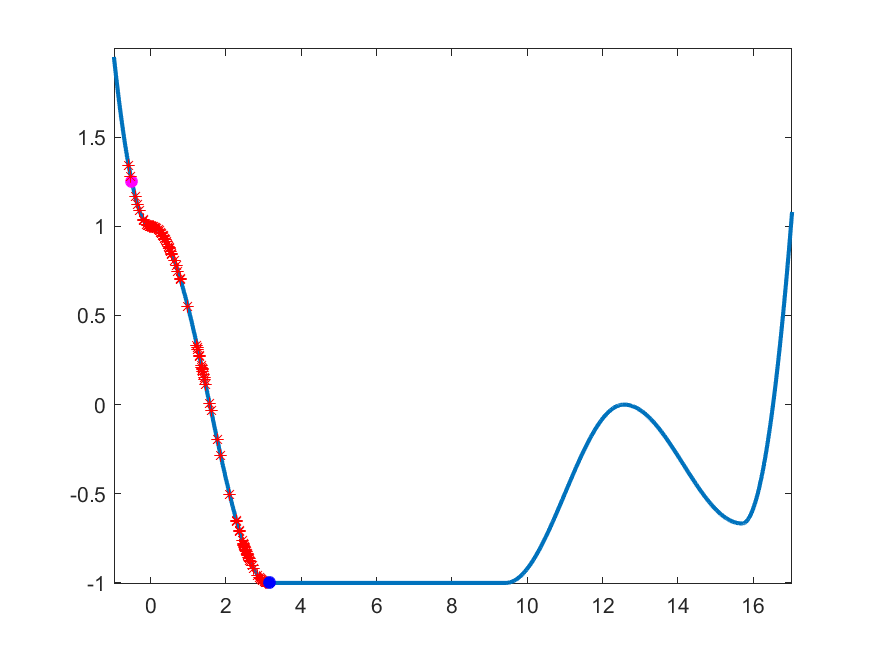}

    \vspace*{-0.5cm}
    
    \caption{Example of a typical situation where process $(\xD_k)_{k\in \Nbb}$ escapes from saddle point $x=0$. Pink point: starting point $x=-1/2$; Blue point: End point belonging to $[\pi, 3\pi]$.}
    \label{fig:escape_saddle_point}
\end{figure}

In our simulations, we choose $b=10$ as a moderate level of multiplicative noise, and add $(c_k)_{k\in \Nbb}$ as an additive one. Despite the higher complexity of the uncertainty model, we obtain slightly better results as the process is able to escape from saddle or local minimizer in all runs (see Figure~\ref{fig:escape_saddle_point}). 


A summary of the behaviour of process $(\xD_k)_{k\in \Nbb}$ under these conditions is reported in Table~\ref{tab2}.

\begin{table}[!t]
    \centering
    \begin{tabular}{||c||c||c|c|}
    \hline 
         $x_0$ & $\sigma$ & Nature of $\xD_{\infty}$ & $F(\xD_k)>\FD_{\infty}~\forall k\in \Nbb$ ? 
        \\ \hline 
        $-1/2$ & $10$  & Global Min & Yes if $F(\xD_k)\neq F(\xD_{\infty})$\\
         & $10^2$  & Global Min & Yes if $F(\xD_k)\neq F(\xD_{\infty})$  \\
         \hline 
          $1$ & $10$  & Global Min & Yes if $F(\xD_k)\neq F(\xD_{\infty})$ \\
          & $10^2$  & Global Min & Yes if $F(\xD_k)\neq F(\xD_{\infty})$  \\
         \hline
          $4\pi + \epsilon $ & $10$  & Global or Local Min &  Yes if $F(\xD_k)\neq F(\xD_{\infty})$ \\
          & $10^2$  & Global or Local Min &  Yes if $F(\xD_k)\neq F(\xD_{\infty})$    \\
         \hline 
    \end{tabular}
    \caption{Asymptotical behaviours of $(\xD_k)_{k\in \Nbb}$ under $\Ebb_k\left[\gD_k^2\right]   \leq 20  F'(\xD_k)^2 + 2\sigma^2 / (k+1)^{(2+2\varepsilon)} ~(k\in \Nbb)$ }
    \label{tab2}
\end{table}

\subsection{Experiments for non-zero $(a_k, b_k, c_k)_{k\in \Nbb}$}

Here we consider the more generic case where none of the sequences $(a_k)_{k\in \Nbb}$, $(b_k)_{k\in \Nbb}$ nor $(c_k)_{k\in \Nbb}$ is equal to $0$.

We suppose that, for every $k\in \Nbb$, the second order moment of $\gD_k$ also possesses a non-zero component relatively to $F(\xD_k)-F_{\min}$. The most difficult assumption to verify is then Ass.~3-\ref{ass:H2_3}, which cannot hold for $(\alpha_k)_{k\in \Nbb}$ and $(b_k)_{k\in \Nbb}$ constant sequences. 
One can show that choosing, for every $k\in \Nbb$, $\alpha_k=\alpha$, $b_k = 3b/(k+1)^2$, and $c_k=3\sigma^2/(k+1)^{2+2\epsilon}$, for $\alpha>0$, $b>0$, $\sigma>0$, and $a_k=2(\alpha^2\beta k^{2+\epsilon})^{-1}$, is sufficient to fulfill Assumptions~\ref{ass:H2}-\ref{ass:H2_1}-\ref{ass:H2_3}. 

We simulate the approximation sequence  $\gD_k =\eD_k^3 \sqrt{F(\xD_k) - F_{\min}}+ \eD_k^1  F'(\xD_k) + \eD_k^2$, where both $\eD_k^1$ and $\eD_k^3$ follows a 1-mean uniform distribution between $[1-\sqrt{3b-1},1+\sqrt{3b-1}]$ and $\left[-\sqrt{3}/\alpha (k+1)^{-(1+\epsilon)}, \sqrt{3}/\alpha (k+1)^{-(1+\epsilon)}\right]$, respectively, and $\eD_k^2\sim\Nc(0, \sigma_k^2)$ with $\sigma_k= \sigma/(k+1)^{(1+\varepsilon)}$.    

A summary of the behaviour of process $(\xD_k)_{k\in \Nbb}$ under these conditions is reported in Table~\ref{tab3}.
The only significant difference of behavior compared with the previous subsection is observed for $(\sigma, x_0)=(10^2, 1)$. In this particular case, numerical errors prevent the algorithm from converging. 

\begin{table}[h]
    \centering
    \begin{tabular}{||c||c||c|c|}
    \hline 
         $x_0$ & $\sigma$ & Nature of $\xD_{\infty}$ & $F(\xD_k)>\FD_{\infty}~\forall k\in \Nbb$ ? 
        \\ \hline 
        $-1/2$ & $10$  & Global Min & Yes if $F(\xD_k)\neq F(\xD_{\infty})$\\
         & $10^2$  & Global Min & Yes if $F(\xD_k)\neq F(\xD_{\infty})$  \\
         \hline 
          $1$ & $10$  & Global Min & Yes if $F(\xD_k)\neq F(\xD_{\infty})$ \\
          & $10^2$  & Non-convergence & Undefined  \\
         \hline
          $4\pi + \epsilon $ & $10$  & Global or Local Min &  Yes if $F(\xD_k)\neq F(\xD_{\infty})$ \\
          & $10^2$  & Global or Local Min &  Yes if $F(\xD_k)\neq F(\xD_{\infty})$    \\
         \hline 
    \end{tabular}
    \caption{Asymptotical behaviours of $(\xD_k)_{k\in \Nbb}$ under $\Ebb_k\left[\gD_k^2\right]   \leq \frac{8}{\alpha^2\beta}(F(\xD_k)-F_{\min})^2+ 40  F'(\xD_k)^2 + 4\sigma^2 / (k+1)^{(2+2\varepsilon)} ~(k\in \Nbb)$.}
    \label{tab3}
\end{table}

\section{Conclusion}
\label{sec5}

In this article we introduce a new theoretical framework to study almost sure convergence of SGD schemes, in a non-convex context.  
We further give numerical illustrations to investigate the behaviour of SGD processes, and the relevancy of the different assumptions necessary to ensure their almost sure convergence. This work is base on the theoretical work we initially conducted in \cite{chouzenoux2023kurdyka}, where we introduced a new KL framework to investigate almost sure convergence of stochastic processes, in a smooth but non-convex context.


\bibliographystyle{abbrv}

\end{document}